\documentclass[oneside,10pt]{article}          
\usepackage[b5paper]{geometry}	    
\usepackage{amsfonts,amsmath,latexsym,amssymb,enumitem} 
\usepackage{amsthm}                
\usepackage{mathrsfs,upref}         
\usepackage{mathptmx}		    
\usepackage{jmi}	            

\newtheorem{theorem}{Theorem}           
               
\newtheorem{corollary}{Corollary}

\theoremstyle{definition}

\newtheorem{rem}{Remark}

\begin{document}

\title{Weighted averages of $\ell^p$ sequences}

\author[L. Bouthat]{Ludovick Bouthat, Javad Mashreghi and Fr\'ed\'eric Morneau-Gu\'erin}

\address{D\'epartement de math\'ematiques et de statistique, Universit\'e Laval, 1045, avenue de la M\'edecine, Qu\'ebec, QC, G1V\;0A6, Canada\\
\email{Ludovick.Bouthat.1@ulaval.ca}}

\address{D\'epartement de math\'ematiques et de statistique, Universit\'e Laval, 1045, avenue de la M\'edecine, Qu\'ebec, QC, G1V\;0A6, Canada\\
\email{Javad.Mashreghi@mat.ulaval.ca}}

\address{D\'epartement \'Education, Universit\'e T\'ELUQ, Qu\'ebec, QC, G1K 9H6, Canada\\
\email{Frederic.Morneau-Guerin@teluq.ca}}

\date{08.06.2023}                               



\begin{abstract}
The objective of the present paper is to establish three Hardy-type inequalities in which the arithmetic mean over a sequence of non-negative real numbers is replaced by some weighted arithmetic mean over some nested subsets of the given sequence of numbers. One of these inequalities stems from a calculation in a paper of Bouthat and Mashreghi on semi-infinite matrices.
\end{abstract}

\maketitle


\section{Introduction}

In his Presidential Address at the meeting of the London Mathematical Society, on November 8, 1928, the great Cambridge mathematician Godfrey Harold Hardy reported the Danish analyst Harald Bohr as saying ``\textit{All analysts spend half their time bunting through the literature for inequalities which they want to use and cannot prove}" \cite[Preface]{MR3676556}.

This assertion -- which, obviously, is a bit of an overstatement -- aptly conveys the idea that inequalities form an essential part of virtually all branches of mathematics. Their usefulness is reflected in the vast literature that exists on the subject.

The Hardy inequality is a prominent example of an inequality that was first derived as a tool to prove a theorem (or, to be more precise, to simplify the proof of some other inequality) but that quickly and unexpectedly took on a life of its own. Indeed, it was not long before the great independent interest of this inequality was realized. It was then studied in its own right by countless mathematicians and acquired several useful variants. It eventually turned out to assume an important role in other areas of mathematics, most notably in the theory of partial differential equations.

The source of motivation for Hardy, when he undertook the research process that led him to state and prove what is now called the Hardy inequality, appears to be a theorem demonstrated by David Hilbert. This theorem, derived in the course of an investigations in the theory of integral equation \cite{MR2256532}, asserts that the double series
$$\sum\limits_{m=1}^\infty \sum\limits_{n=1}^\infty \frac{a_m a_n}{m+n}$$
with positive terms $a_m \geq 0$ converges whenever $\sum_{m=1}^\infty a_m^2$ is convergent. More precisely, Hilbert showed that 
\begin{equation}\label{EqHilbert}\sum\limits_{m=1}^\infty \sum\limits_{n=1}^\infty \frac{a_m a_n}{m+n} \,\leq\, 2\pi \sum\limits_{m=1}^\infty a_m^2. \end{equation}

At least four essentially different proofs of this beautiful fact were published in a few years period :
\begin{enumerate}
\item[(1)] The original proof of Hilbert, which depends upon the theory of Fourier's series, was published in \cite{MR0056184} in 1906. It was also outlined by his doctoral student H. Weyl in 1908 in his doctoral dissertation \cite{WeylThesis} as well as in \cite{MR1511502}.
\item[(2)] An elementary yet unnatural proof was obtained by F. W. Wiener in 1910 and it appeared in \cite{MR1511566}.
\item[(3)] I. Schur provided a proof depending upon the theory of quadratic and bilinear forms in an infinity of variables in 1911 in \cite{MR1580823}.
\item[(4)] In the same paper, Schur also derived a second proof relying on a change of variables in a double integral. By Hardy’s standard, it is this latter proof that was ``the most elegant of all”  \cite{MR1544414}.
\end{enumerate}

The determination that the constant $2\pi$ in \eqref{EqHilbert} may be replaced by $\pi$, in which case the inequality is sharp, is also due to Schur.
However, none of the aforementioned proofs being considered ``\textit{as simple and elementary as might be desired}” by Hardy. He then worked yet another proof in 1918 which seemed to him ``\textit{to lack nothing in simplicity}” \cite{MR1544414}.

Hardy's elementary proof, which was published the following year in \cite{hardy1919notes}, consists in deriving Hilbert’s inequality as an easy corollary to the following statement: if $\sum_{n=1}^\infty a_n^2$ is convergent, then so is 
$$\sum\limits_{n=1}^\infty \left(\frac{a_1 + a_2 + \cdots + a_n}{n}\right)^{\!\!2}.$$

Hardy perspicaciously reflected in \cite{MR1544414} that this  inequality (which he had already  hinted at under a slightly different form and with only a sketch of proof in a 1915 paper \cite{hardy1915notes}) ``\textit{seems to be of some interest in itself}”.

Upon hearing about the above inequality from Hardy, a frequent collaborator, Marcel Riesz came up with an alternative proof leading to the following generalization : \textit{if $p > 1$, $a_n \geq 0$ and $\sum_{n=1}^\infty a_n^p$ is convergent, then}
\begin{equation}\label{EqHardyRiesz}\sum\limits_{n=1}^\infty\left(\frac{a_1 + a_2 + \cdots + a_n}{n}\right)^{\!p} \leq\,  \left( \frac{p^2}{p-1}\right)^{\!p} \sum\limits_{n=1}^\infty a_n^p.\end{equation}

Riesz’s proof as well as a short historical account of the events that led to its discovery was published in 1920 by Hardy \cite{MR1544414}. Interestingly, towards the end of the text Hardy indicates that he is well aware that the inequality \eqref{EqHardyRiesz} remains true for all relevant sequences if the constant $\left( p^2/(p-1)\right)^p$ is replaced by some smaller constant. Somewhat surprisingly, although everything indicates that he believed that the smaller constant $\left( p/(p-1)\right)^p$ was sharp (one reason for this was a proof “sent to him by Prof. Schur by letter” that this was at least the case for $p=2$), Hardy stopped short of asserting it outright. Instead, he goes on saying ``\textit{It would be interesting to know what is the best constant which can appear on the right hand side} [\textit{of} \eqref{EqHardyRiesz}]”. This shortcoming was remedied when E. Landau provided Hardy with a proof that the constant $\left( p/(p-1)\right)^p$ is indeed sharp in a letter to Hardy -- dated June 21, 1921 -- that was published five years later \cite{MR1575105}.

In his 1920 note, also provided (without proof) was the following integral formulation : \textit{if $p > 1$, $a > 0$, $f(x) \geq 0$, and $\int_a^\infty f(x)^p \; dx$ is convergent, then}
$$ \int_a^\infty \!\left( \frac{\int_a^x f(t)\; dt}{x}\right)^{\!p}  dx ~\leq~ \left(\frac{p}{p-1}\right)^{\!p} \!\int_a^\infty \!f^p(x)\; dx.$$

A plethora of alternative proofs of Hardy’s inequality have been given by various authors in the following years. See for example the work of T. A. A. Broadbent \cite{MR1574000}, E. B. Elliott \cite{MR1574962}, K. Grandjot \cite{MR1574375}, G. H. Hardy himself \cite{hardy1925notes}, T. Kaluza and G. Szeg\"{o} \cite{MR1574830}, and K. Knopp \cite{MR1574165}.

Nowadays, the literature on the Hardy inequality, its variations, its generalizations, and its applications is so vast (on an indicative basis, a search on MathSciNet for papers containing the words ``Hardy inequality" in the title returned more than a thousand results to display) that an exhaustive list of reference is simply impossible.

Let it therefore suffice to cite a few prominent contributions, such as the following monographs \cite{MR3408787, MR1747888, MR2091115, MR3287251,  MR1069756}.

For a simple and accessible presentation of the Hardy inequality as well as some of its extensions and applications, we refer to \cite{masmoudi2011hardy}. For a thorough overview of Hardy inequalities on time scales, see \cite{MR3561405}. For a comprehensive look at weighted Hardy inequalities, we point towards \cite{MR3676556}. As for an extensive treatment of Hardy inequalities and closely related topics from the point of view of Folland and Stein's homogeneous groups, see \cite{ruzhansky2019hardy}. For a very detailed presentation of the history of the Hardy inequality, see \cite{MR2256532, MR2351524}. Lastly, the following article offers a quite enlightening uchronia \cite{MR2925639}.

The objective of the present paper is to establish three Hardy-type inequalities in which the arithmetic mean is replaced by some weighted arithmetic mean. The form of these inequalities closely resemble that of the classical Hardy inequality. But, strictly speaking, none of them constitute a proper generalization of \eqref{EqHardyRiesz} as they do not entail this equation as a special case.

Following a brief overview, in Section \ref{ContexSection}, of the context surrounding our results,  we present in Section \ref{MainSection} the two results that form the core of this article. In Section \ref{MainProof}, we provide a proof of our first main theorem. The following Section contains a complete proof of a corollary applying to lacunary series. Finally, in Section \ref{SecondThmProof} we give a proof of our second main theorem.


\section{Some elements of context}\label{ContexSection}

In studying semi-infinite matrices which represent a linear bounded operator on the sequence space $\ell^2$, the inequality
\begin{equation}\label{E:4n-ineq}
\sum_{n=1}^{\infty} \frac{1}{4^n} \left| \sum_{j=1}^{4^n} a_j \right|^2
\leq\,
C \sum_{n=1}^{\infty} |a_n|^2,
\end{equation}
where $C$ is a positive constant, appeared in the calculation \cite[Example 1]{MRtmp}. The verification of this simple-looking inequality was not straightforward. However, after completing the first steps then applying the inequality in the initial context, we naturally asked whether the inequality  \eqref{E:4n-ineq} is a special case of a more general one or not. The purpose of the present note is to establish a generalization of \eqref{E:4n-ineq}.

To set the stage, let us introduce some notation: let $(a_n)_{n \geq 1}$ be a sequence of complex numbers in $\ell^2$, i.e.,
\[
\sum_{n=1}^{\infty} |a_n|^2 \,<\, \infty.
\]
Given a strictly increasing sequence of positive numbers $(n_k)_{k \geq 1}$, let $m_k \geq 0$, $M_k>0$, $k \geq 1$, and consider the average-like combinations
\[
b_k \,:=\, \frac{1}{M_k} \sum_{j=1}^{n_k} m_j a_j, \qquad (k \geq 1).
\]
We seek conditions on $m_k$, $M_k$, and $n_k$, such that the sequence $(b_n)_{n \geq 1}$ belongs to $\ell^2$ and that its norm is controlled by the norm of $(a_n)_{n \geq 1}$, i.e.,
\[
\sum_{n=1}^{\infty} |b_n|^2 \,\leq\, C \sum_{n=1}^{\infty} |a_n|^2,
\]
or, more explicitly,
\begin{equation}\label{E:lacun-1}
\left| \frac{1}{M_1} \sum_{j=1}^{n_1} m_j a_j \right|^2 \!\!+\, \left| \frac{1}{M_2} \sum_{j=1}^{n_2} m_j a_j \right|^2 \!\!+\, \cdots \,\leq\, C \sum_{n=1}^{\infty} |a_n|^2.
\end{equation}
Clearly, the inequality \eqref{E:4n-ineq} is a special case of \eqref{E:lacun-1}.


\section{Main results}\label{MainSection}

In this section, we present a variation of \eqref{E:lacun-1}, even slightly more general, in which the parameters $m_j$, $M_j$ and $n_j$ are intertwined. Let $\mathbb{N}$ denote the set of positive integers, and let
\[
\mathbb{N} ~=~ N_1 \cup N_2 \cup \cdots
\]
be a partition of $\mathbb{N}$. Put
\[
\mathbf{N}_n \,:=\, N_1\cup\cdots\cup N_n, \qquad (n \geq 1).
\]
The $\mathbf{N}_n$'s form an increasing nested sequence of sets that exhaust $\mathbb{N}$, i.e., 
\[
\mathbf{N}_1 \,\varsubsetneqq\, \mathbf{N}_2 \,\varsubsetneqq\, \mathbf{N}_3 \,\varsubsetneqq\, \cdots,
\]
and
\[
\mathbb{N} \,=\, \bigcup\limits_{n=1}^\infty \mathbf{N}_n.
\]
Let $(m_n)_{n \geq 1}$ be a sequence of weights (positive numbers) and $p > 1$ and $q >1$ be conjugate exponents, i.e., $\frac{1}{p} + \frac{1}{q}=1$, and let us define
\begin{equation}\label{E:def-wj}
w_n \,:=\, \left(\sum_{j \in N_n} m_j^q \right)^{\!\!1/q}\!\!,\qquad (n \geq 1),
\end{equation}
and
\begin{equation}\label{E:2}
M_n \,:=\, \sum_{j=1}^{n} w_j, \qquad (n \geq 1).
\end{equation}

In the following, we will consider the sequence
\[
\frac{1}{M_n}\sum_{j \in \mathbf{N}_n} m_ja_j,
\]
which can be interpreted as some sort of weighted arithmetic average of $(a_n)_{n \geq 1}$. However, note that it is in fact not a usual weighted arithmetic average since in general, $\sum_{j \in \mathbf{N}_n} m_j \neq M_n$.

As we shall see, the quantity that connects the above parameters and that plays a major role below is
\begin{equation}\label{E:3}
\rho \,:=\, \sup_{n \geq 1} \, \left(w_n \sum_{j=n}^{\infty}\frac{1}{M_j}\right)^{\!\!1/p}.
\end{equation}

\begin{theorem} \label{T:main-ineq}
Let $(a_n)_{n \geq 1}$ be a sequence of complex numbers, and let $(m_n)_{n \geq 1}$ be a sequence of weights. Define $M_n$ and $\rho$ as in \eqref{E:2} and \eqref{E:3} respectively, and assume that $\rho$ is finite.
Then
\begin{align}\label{eq:main}
	\left(\sum_{n=1}^{\infty} \left|\frac{1}{M_n}\sum_{j \in \mathbf{N}_n} m_ja_j \right|^p\right)^{\!\!1/p}
	\!\!\leq~ \rho  \left(\sum_{n=1}^{\infty} |a_n|^p\right)^{\!\!1/p}.
\end{align}
\end{theorem}

\begin{rem}\label{R:1}
	More generally, if $(M_n)_{n\geq1}$ is any given sequence of positive numbers (as opposed to being defined by \eqref{E:2}), then \eqref{eq:main} remain true if $\rho$ is replaced by the constant
	\[
	\sup_{n \geq 1} \, \left(w_n \sum_{j=n}^{\infty}\frac{\left( w_1 + \cdots +w_k \right)^{p-1}}{M_{k}^{p}}\right)^{\!\!1/p}\!\!,
	\]
	whenever it is finite. We shall prove this fact and Theorem~\ref{T:main-ineq} in one stroke in the following section. 
	
	In this case, observe that for $m_j=1$ ($j \geq 1$), $M_n=n^{1+\varepsilon}$ ($\varepsilon>0$) and $N_n = \{n\}$ ($n \geq 1$), we obtain an Hardy-like inequality, namely an inequality of the same form as the classical Hardy inequality \eqref{EqHardyRiesz} but with the denominator on the left-hand side replaced by $n^{1+\varepsilon}$.
\end{rem}

We say that a sequence $(n_k)_{k \geq 1}$ of positive integers satisfy the Hadamard gap condition if
\[
\frac{n_{k+1}}{n_{k}} \geq r > 1, \qquad (k \geq 1).
\]
Whenever this is the case, $(n_k)_{k \geq 1}$ is called a lacunary sequence of ratio $r$ (see \cite{ hadamard1892essai, MR1738432}).

\begin{corollary} \label{C:lacun-ineq}
	Let $(a_n)_{n \geq 1}$ be a sequence of complex numbers, and let $(n_k)_{k \geq 1}$ be a lacunary series with ratio $r$. If $p,q\in(1,\infty)$ are such that $\frac{1}{p}+\frac{1}{q}=1$, then
	\[
	\left( \sum_{k=1}^{\infty} \left| \frac{1}{n_k^{1/q}} \sum_{j=1}^{n_k} a_j \right|^p\right)^{\!\!1/p}
	\!\!\leq~
	 \frac{r^{1/q}}{r^{1/q}-1}  \left(\sum_{n=1}^{\infty} |a_n|^p\right)^{\!\!1/p}.
	\]
\end{corollary}

If the sequence $n_k$ is geometric, i.e., $n_k = b^k$ for some $b>1$, then in the particular case of $p=2$ we can go even further and establish the following result in which the constant is optimal.

\begin{theorem} \label{T:geometric-seq}
	Let $(a_n)_{n \geq 1}$ be a sequence of complex numbers and let $b>1$ be an integer. Then
	\begin{equation}\label{optimal_eq}
		\left( \sum_{k=1}^{\infty} \frac{1}{b^k}  \left| \sum_{j=1}^{b^k} a_j \right|^2\right)^{\!\!\!1/2}
		\!\!\leq~
		\frac{\sqrt{b}+1}{\sqrt{b-1}} \left(\sum_{n=1}^{\infty} |a_n|^2\right)^{\!\!1/2}.
	\end{equation}
	Moreover, the constant $\frac{\sqrt{b}+1}{\sqrt{b-1}}$ is sharp and the above inequality is strict, except if $(a_n)_{n\geq1}$ is the null sequence.
\end{theorem}


\section{Proof of Theorem \ref{T:main-ineq}}\label{MainProof}

	Observe that
	\begin{equation*}
		\begin{aligned}
			\sum_{j \in \mathbf{N}_n} m_ja_j ~&= \sum_{j \in N_1\cup\cdots\cup N_n} \!\! m_ja_j\\
			&=~ \sum_{j \in N_1} m_ja_j+
			\sum_{j \in N_2} m_ja_j+\cdots+
			\sum_{j \in N_n} m_ja_j.
		\end{aligned}
	\end{equation*}
	Applying Hölder's inequality individually to each of the latter sums yields
	\begin{eqnarray*}
		& & \sum_{j \in \mathbf{N}_n} m_ja_j\\
	&\leq& \left(\sum_{j \in N_1} m_j^q\right)^{\!\!1/q} \!\!\!\!\!\cdot \left(\sum_{j \in N_1} a_j^p\right)^{\!\!1/p}
	\!\!\!\!+\cdots+
	\left(\sum_{j \in N_n} m_j^q\right)^{\!\!1/q} \!\!\!\!\!\cdot \left(\sum_{j \in N_n} a_j^p\right)^{\!\!1/p} \\
	&=& w_1 \left(\sum_{j \in N_1} a_j^p\right)^{\!\!1/p}
	\!\!\!\!+\cdots+
	w_n \left(\sum_{j \in N_n} a_j^p\right)^{\!\!1/p}\\
	&=& \sum\limits_{k=1}^n w_k \cdot \left(\sum_{j \in N_k} a_j^p\right)^{\!\!1/p}.
	\end{eqnarray*}
	The next non-trivial step consists of reformulating the latter term as
	\[
	\sum\limits_{k=1}^n w_k^{1/q} \cdot \left[ w_k^{1/p} \left(\sum_{j \in N_k} a_j^p\right)^{\!\!1/p} \right],
	\]
	and to then apply, once more, Hölder's inequality to obtain
	\begin{equation}
		\begin{aligned}\label{E:sum-ineq-average}
		\sum_{j \in \mathbf{N}_n} m_ja_j
		\,\leq\, \left( w_1 + \cdots +w_n \right)^{1/q} \left( w_1 \sum_{j \in N_1} a_j^p + \cdots + w_n \sum_{j \in N_n} a_j^p\right)^{\!\!1/p}\!\!\!.
		\end{aligned}
	\end{equation}

	By \eqref{E:3}, we can rewrite this estimation as
	\[
	\left(\frac{1}{M_n}\sum_{j \in \mathbf{N}_n} m_ja_j \right)^{\!\!p} 
	\leq\, W_{1,n} \sum_{j \in N_1} a_j^p + \cdots + W_{n,n} \sum_{j \in N_n} a_j^p,
	\]
	where
	\[
	W_{n,k} \,:=\, \frac{w_n \left( w_1 + \cdots +w_k \right)^{p-1}}{M_{k}^{p}}.
	\]

	If we write these inequalities for $n\geq 1$ and add them up, we see that
	\begin{eqnarray*}
		\sum_{n=1}^{\infty} \left(\frac{1}{M_n}\sum_{j \in \mathbf{N}_n} m_ja_j \right)^{\!\!p}
		&\leq& \sum_{n=1}^{\infty} \left( W_{1,n} \sum_{j \in N_1} a_j^p + \cdots + W_{n,n} \sum_{j \in N_n} a_j^p\right) \\
		&=& \sum_{n=1}^{\infty} \left[ \left( \sum_{k=n}^{\infty} W_{n,k} \right)  \sum_{j \in N_n} a_j^p\right]\\
		&\leq& \sum_{n=1}^{\infty} \rho^p \sum_{j \in N_n} a_j^p\\
		&=& \rho^p  \sum_{n=1}^{\infty} a_n^p
	\end{eqnarray*}
	where $\rho=\sup_{n \geq 1} \, \left(\sum_{j=n}^{\infty} W_{n,j}\right)^{\!1/p}$. Since by \eqref{E:2}, $M_n=w_1+\cdots+w_n$, we have $W_{n,j} = \frac{w_n}{M_j}$ and we are done. $\hfill\qed$


\section{Proof of Corollary \ref{C:lacun-ineq}}\label{SectionCorollary}

Let $N_1=\{1,2,\dots,n_1\}$ and
	\[
	N_k \,=\, \{n_{k-1}+1,n_{k-1}+2,\dots,n_k\}, \qquad (k \geq 2),
	\]
	$m_k=1$ for all $k \geq 1$, and $M_k=n_k^{1/q}$. Then, by \eqref{E:def-wj},
	\[
	w_1 \,=\, \left(\sum_{j\in N_1} m_j^q \right)^{\!\!1/q} \!\!\!=\, \left(\sum_{j=1}^{n_1} m_j^q \right)^{\!\!1/q} \!\!\!=\, n_{1}^{1/q}
	\]
	and
	\[
	w_k \,=\, \left(\sum_{j\in N_k} m_j^q \right)^{\!\!1/q} \!\!\!=\, \left(\sum_{j=n_{k-1}+1}^{n_k} m_j^q \right)^{\!\!1/q} \!\!\!=\, (n_k-n_{k-1})^{1/q},  \qquad (k \geq 2).
	\]
	In any case,
	\begin{equation}\label{E:wk-nk}
		w_k \,\leq\, n_k^{1/q}, \qquad  (k \geq 1)
	\end{equation}
	By Remark \ref{R:1}, we thus have
	\[
	\left( \sum_{k=1}^{\infty} \left| \frac{1}{n_k^{1/q}} \sum_{j=1}^{n_k} a_j \right|^p\right)^{\!\!1/p} \!\!\leq\, \rho \left(\sum_{k=1}^{\infty} |a_k|^p\right)^{\!\!1/p}
	\]
	where  $\rho=\sup_{m \geq 1} \, \left(\sum_{j=m}^{\infty} \frac{w_m \left( w_1 + \cdots +w_j \right)^{p-1}}{M_{j}^{p}}\right)^{\!1/p}$. However, by \eqref{E:wk-nk} and using the fact that $n_k \leq r^{k-j}n_j$, we find that
	\begin{align*}
		\frac{w_m \left( w_1 + \cdots +w_j \right)^{p-1}}{M_{j}^{p}} \,&\leq\, \frac{n_m^{1/q} \left( r^{(1-j)/q}n_j^{1/q} + \cdots +n_j^{1/q} \right)^{p-1}}{n_{j}^{p/q}} \\
		&=\, \left(\frac{n_m}{n_j} \right)^{\!\!1/q}\!\! \left( 1 + r^{-1/q} + \cdots +r^{-(j-1)/q} \right)^{p-1} \\
		&\leq\, r^{(m-j)/q} \left( \frac{r^{1/q}}{r^{1/q}-1} \right)^{p-1}.
	\end{align*}
	Hence, 
	\begin{align*}
		\sum_{j=m}^{\infty} \frac{w_m \left( w_1 + \cdots +w_j \right)^{p-1}}{M_{j}^{p}} \,&\leq\, \left( \frac{r^{1/q}}{r^{1/q}-1} \right)^{\!p-1}  \sum_{j=m}^{\infty} r^{(m-j)/q} \\
		&=\, \left( \frac{r^{1/q}}{r^{1/q}-1} \right)^{\!p-1}  \sum_{s=0}^{\infty} r^{-s/q} \\
		&=\, \left( \frac{r^{1/q}}{r^{1/q}-1} \right)^{\!p-1} \!\! \left(\frac{r^{1/q}}{r^{1/q}-1}  \right) \\
		&=\, \left( \frac{r^{1/q}}{r^{1/q}-1} \right)^{\!p}
	\end{align*}
	and it follows that $\rho \leq \frac{r^{1/q}}{r^{1/q}-1}$. $\hfill\qed$


\section{Proof of Theorem \ref{T:geometric-seq}}\label{SecondThmProof}

	The notations are as in the Corollary \ref{C:lacun-ineq} and the proof has some overlap. However, we need more precise calculation here.
	By \eqref{E:def-wj},
	\[
	w_1 \,=\, \left(\sum_{j\in N_1} m_j^2 \right)^{\!\!1/2} \!\!\!=\, \left(\sum_{j=1}^{n_1} m_j^2 \right)^{\!\!1/2} \!\!\!=\, \sqrt{n_1} \,=\, \sqrt{b}
	\]
	and, for $k \geq 2$,
	\begin{align*}
		w_k \,&=\, \left(\sum_{j\in N_k} m_j^2 \right)^{\!\!1/2} \!\!\!=\, \left(\sum_{j=n_{k-1}+1}^{n_k} m_j^2 \right)^{\!\!1/2} \!\!\!=\, \sqrt{n_k-n_{k-1}} \\
		&=\, \sqrt{b^k-b^{k-1}} \,=\, \sqrt{b-1} \, \sqrt{b}^{k-1}.
	\end{align*}
	Hence,
	\begin{align*}
		w_1+\cdots + w_k \,&=\, \sqrt{b} + \sqrt{b-1} (\sqrt{b} + \sqrt{b}^2 + \cdots + \sqrt{b}^{k-1}) \\
		&=\, \sqrt{b} + \sqrt{b-1} \frac{\sqrt{b}^k-\sqrt{b}}{\sqrt{b}-1}.
	\end{align*}
	Therefore, by \eqref{E:sum-ineq-average},
	\begin{eqnarray*}
		\left|\sum_{j =1}^{n_k} a_j \right|^{2}
		&\!\!\leq&\!\left( w_1 + \cdots +w_n \right) \left( w_1 \sum_{j \in N_1} a_j^2 + \cdots + w_n \sum_{j \in N_n} a_j^2\right) \\
		&\!\!=&\! \bigg( \sqrt{b} + \sqrt{b-1} \frac{\sqrt{b}^k-\sqrt{b}}{\sqrt{b}-1} \bigg) \\
		&&~~~\! \cdot\left( \sqrt{b} \sum_{j \in N_1} a_j^2 + \cdots + \sqrt{b-1} \sqrt{b}^{k-1} \sum_{j \in N_k} a_j^2\right),
	\end{eqnarray*}
	which gives
	\[
	\left| \frac{1}{\sqrt{n_k}}\sum_{j=1}^{n_k} a_j \right|^2 \!\leq \frac{\sqrt{b} + \sqrt{b-1} \frac{\sqrt{b}^k-\sqrt{b}}{\sqrt{b}-1}}{b^k} \left( w_1 \sum_{j \in N_1} a_j^2 + \cdots + w_k \sum_{j \in N_k} a_j^2\right).
	\]

	For $k \geq 1$, write
	\[
	W_k \,:=\, \frac{b^k}{\sqrt{b} + \sqrt{b-1} \frac{\sqrt{b}^k-\sqrt{b}}{\sqrt{b}-1} }.
	\]
	Then,
	\begin{equation*}
	\begin{aligned}	
	\sum_{k=1}^{\infty}  \left| \frac{1}{\sqrt{n_k}}\sum_{j=1}^{n_k} a_j \right|^2 &\leq\, \sum_{k=1}^{\infty} \frac{1}{W_k}\left( w_1 \sum_{j \in N_1} a_j^2 + \cdots + w_k \sum_{j \in N_k} a_j^2\right) \\&=\, \sum_{k=1}^{\infty} \left[ \left(w_k \sum\limits_{j=k}^\infty \frac{1}{W_j}\right) \sum\limits_{j \in N_k}a_j^2 \right] \\&\leq\, \rho' \sum\limits_{k=1}^\infty \sum\limits_{j \in N_k} a_j^2  \,=\, \rho' \sum_{j =1}^{\infty} a_j^2,
	\end{aligned}
	\end{equation*}
	where $\rho'= \sup\limits_{k\geq 1} \left[ w_k \sum_{j=k}^{\infty}\frac{1}{W_j} \right]$. But
	\begin{align*}
		w_1 \sum_{j=1}^{\infty}\frac{1}{W_j} \,&=\, \sqrt{b} \sum_{j=1}^{\infty}\frac{1}{W_j} \,=\, \sqrt{b} \sum_{j=1}^{\infty} \frac{ \sqrt{b} + \sqrt{b-1} \frac{\sqrt{b}^k-\sqrt{b}}{\sqrt{b}-1}}{b^k} \\
		&=\, \frac{\sqrt{b-1}+\sqrt{b}(\sqrt{b}-1)}{b-1} \,<\, \frac{(1+\sqrt{b})^{2}}{b-1}
	\end{align*}
	and, for $k \geq 2$,
	\begin{align*}
		w_k \sum_{j=k}^{\infty}\frac{1}{W_j} \,&=\, \sqrt{b-1}\sqrt{b}^{k-1} \sum_{j=k}^{\infty} \frac{ \sqrt{b} + \sqrt{b-1} \frac{\sqrt{b}^k-\sqrt{b}}{\sqrt{b}-1}}{b^k} \\
		&=\, \frac{(1+\sqrt{b})^2+\sqrt{b}^{2-k}(\sqrt{b-1}-\sqrt{b}-1)}{b-1} \\
		&<\, \frac{(1+\sqrt{b})^2}{b-1}.
	\end{align*}

	Hence
	\[
	\rho' \,<\, \frac{(1+\sqrt{b})^{2}}{b-1},
	\]
	thus proving the required inequality. Moreover, since the above inequality is strict, this also prove the fact that \eqref{optimal_eq} is always strict except if $(a_n)_{n\geq 1}$ is the null sequence.

	To prove that this constant is sharp, we provide a specific non-trivial example for which equality is attained. Let
	\[
	a_j \,=\, r^{-k}, \qquad (j\in N_k, \, k \geq 1),
	\]
	where $r>\sqrt{b}$ (a restriction that will be needed later). Then
	\begin{align*}
		\frac{1}{\sqrt{n_k}}\sum_{j=1}^{n_k} a_j  \,&=\, \frac{1}{\sqrt{b}^k} \left( b \cdot r^{-1} + \sum_{j=2}^k (b^j-b^{j-1}) \cdot r^{-j}  \right)\\
		&=\, \frac{1}{\sqrt{b}^k} \left( \frac{b}{r} + \frac{b-1}{b}\sum_{j=2}^k \left(\frac{b}{r}\right)^{j} \right) \\
		&=\, \frac{1}{\sqrt{b}^k} \left( \frac{b}{r} + \frac{b-1}{b}\cdot \frac{b\left(b-r(b/r)^k\right)}{r(r-b)} \right) \\
		&=\, \frac{b}{r\sqrt{b}^k} \left(  \frac{r-1-(b-1)(b/r)^{k-1}}{r-b} \right).
	\end{align*}
	Hence,
	\begin{equation*}
		\begin{aligned}
			&&&\left| \frac{1}{\sqrt{n_k}}\sum_{j=1}^{n_k} a_j \right|^2 \\
			&=&&\frac{b^2}{r^2 b^k} \left(  \frac{r-1-(b-1)(b/r)^{k-1}}{r-b} \right)^{\!2} \\
			&=&&\frac{b^2}{r^2 b^k} \cdot \frac{(r-1)^2-2(r-1)(b-1)(b/r)^{k-1}+(b-1)^2 (b/r)^{2k-2}}{(r-b)^2}\\
			&=&&\frac{(r-1)^2}{(r-b)^2}  \bigg( \frac{b^2}{r^2}\! \left(\frac{1}{b}\right)^{\!\!k} \!- 2\frac{b(b-1)}{r(r-1)}\! \left(\frac{1}{r}\right)^{\!\!k} \!+\frac{(b-1)^2}{(r-1)^2}\! \left(\frac{b}{r^2}\right)^{\!\!k} \bigg).
		\end{aligned}
	\end{equation*}
	Summing up these expressions with respect to $k\geq 1$ gives
	\[
	\sum_{k=1}^{\infty}  \left| \frac{1}{\sqrt{n_k}}\sum_{j=1}^{n_k} a_j \right|^2 \!=\, \frac{b^2}{r^2 (r-b)^2} \left( \frac{(r-1)^2}{b-1} - \frac{2r(b-1)}{b} + \frac{r^2(b-1)^2}{b(r^2-b)} \right).
	\]
	Moreover,
	\begin{align*}
		\sum_{n=1}^{\infty} |a_n|^2 \,&=\, b\cdot r^{-2} +  \sum_{n=2}^{\infty} (b^n-b^{n-1})\cdot r^{-2n} \\
		&=\, \frac{b}{r^2} +  \frac{b-1}{b}\sum_{n=2}^{\infty} \left(\frac{b}{r^2}\right)^{n} \\
		&=\, \frac{b}{r^2} +  \frac{b-1}{b}\cdot \frac{b^2}{r^2(r^2-b)} \\
		&=\,
		\frac{b (r^2 - 1)}{r^2 (r^2 - b)}.
	\end{align*}
	Hence, in the light of our inequality, 
	\[
	C^2 \,\geq\, \frac{\sum_{k=1}^{\infty}  \left| \frac{1}{\sqrt{n_k}}\sum_{j=1}^{n_k} a_j \right|^2}{\sum_{n=1}^{\infty} |a_n|^2},
	\]
	and we have
	\begin{align*}
		C^2 \,&\geq\, \frac{b(r^2-b)}{(r^2-1) (r-b)^2} \left( \frac{(r-1)^2}{b-1} - \frac{2r(b-1)}{b} + \frac{r^2(b-1)^2}{b(r^2-b)} \right) \\
		&=\, \frac{b(r^2-b)}{(r^2-1) (r-b)^2} \left( \frac{(r-1)^2}{b-1} - \frac{2r(b-1)}{b} \right) + \frac{r^2(b-1)^2}{(r^2-1)(r-b)^2}.
	\end{align*}
	It suffices to let $r\to \sqrt{b}$ to deduce
	\begin{equation*}
		C ~\geq~ \frac{\sqrt{b}+1}{\sqrt{b-1}},
	\end{equation*}
	and we are done. $\hfill\qed$



\bibliographystyle{plain}
\bibliography{CRreferences}

\end{document}